# Polygonal billiards and "optical tori"

Eduardo Díaz-Miguel

I.E.S. Salvador Rueda, Departamento de Matemáticas, 29006, Málaga, Spain.

Departamento de Física Aplicada I, Facultad de Ciencias, Universidad de Málaga, 29071-Málaga, Spain.

## Abstract

We give an optical physicist view of the problem of the trajectories in a polygonal billiard using only basic facts of Optics and the theory of functions of a complex variable. This approach allows us to stablish a certain correspondence between n-gon billiards and one-holed 2n-punctured flat tori. Therefore the existence of periodic trajectories in a certain polygon becomes the problem of the existence of closed geodesics in its associated torus.

## 1. An optical problem

Let a light ray propagate in a certain domain of a homogeneous and isotropic bidimensional medium. This means that its refractive index n(x,y) is constant and the velocity of the light rays does not depend on the direction of propagation. As an example, let us suppose that the boundary C of that domain T is a square whose sides are mirrors (a quadrilateral billiard). In Fig.(1), the trajectory of the ray is the polygonal 1-2-3-4-5-6-. Now we apply to T a geometrical transformation, given by a holomorphic function of a complex variable, which transforms T into $T^*$ and the boundary C of T into the boundary $C^*$ of $T^*$. Let's identify the points of $T \cup C$ and $T^* \cup C^*$ with the complex numbers z = x+iy and w = u+iv, respectively; so that w = f(z) = u(x,y)+iv(x,y). As f is holomorphic it is conformal, therefore the transformed trajectory $1^*$-$2^*$-$3^*$-$4^*$-$5^*$-$6^*$- also satisfy the law of reflection in $T^* \cup C^*$. For example, Fig. (2) is the transformed of Fig. (1) by the mapping:

$$w = f(z) = i e^{-z} \quad \Rightarrow \quad \begin{cases} u(x,y) = e^{-x} \sin y \\ v(x,y) = e^{-x} \cos y \end{cases} \quad (1)$$

Obviously, the curvilinear polygonal $1^*$-$2^*$-$3^*$-$4^*$-$5^*$-$6^*$- cannot be the trajectory of a light ray in any homogeneous medium in $T^*$. However we can state the following

## Problem

Is it possible to fill the domain $T^*$ with a non-homogeneous medium, with refractive index n(u,v), such that the curvilinear polygonal trajectories of the light rays in this medium coincide with the transformed of the rectilinear polygonal trajectories in T under the action of f?

Let us see that the answer is affirmative. Accordingly, such a refractive index could be called a *flat index*.

## 2. The solution

Fermat´s principle asserts [1, 2] that the light trajectories in a medium whose refractive index is n(u,v) are the geodesics of a Riemannian space whose (conformal) metric is given by the line element:

$$ds^2 = n(u,v)^2(du^2 + dv^2) \qquad (2)$$

In T the straight lines are the geodesics of the Euclidean metric $ds^2 = dx^2 + dy^2$.
A simple calculation, that takes into account the Cauchy-Riemann equations, shows that this metric in the (u,v) coordinates is

$$ds^2 = \left(\left(\frac{\partial x}{\partial u}\right)^2 + \left(\frac{\partial x}{\partial v}\right)^2\right)(du^2 + dv^2) \qquad (3)$$

But the first parenthesis of (3) is the square modulus of the derivative of the inverse function:

$$\frac{d f^{-1}(w)}{dw} = \frac{dz}{dw} = \frac{\partial x}{\partial u} + i\frac{\partial x}{\partial v}$$

So the searched refractive index is related to the mapping f trough the relationship:

$$n(u,v) = \left|\frac{dz}{dw}\right| = \left|\frac{dw}{dz}\right|^{-1} \qquad (4)$$

The refractive index which must fill the domain $T^*$ of Fig. (2) is, according to (1) and (4):

$$n(u,v) = \left|\frac{d}{dw}(-i\ln w)\right| = \left|\frac{1}{w}\right| = \frac{1}{\sqrt{u^2+v^2}} = \frac{1}{r^*}, \qquad (5)$$

where $r^*$ is the distance of the point (u,v) to the origin of coordinates. Let us observe that the square T of Fig. (1) has been chosen so that $T^*$ does not contain the point (0,0). Then n(u,v) is not singular. Moreover $r^* \leq 1$ in $T^*$ implies $n(u,v) \geq 1$.

   It is immediate to prove that a flat index given by Equation (4) is the exponential of a harmonic function:

$$\Delta(\ln n) = 0 \qquad (6)$$

 This is a consequence of the fact that the metric (3) is the Euclidean one expressed in the (u,v) coordinates; so its Gaussian curvature is cero and Equation (6) easily follows.

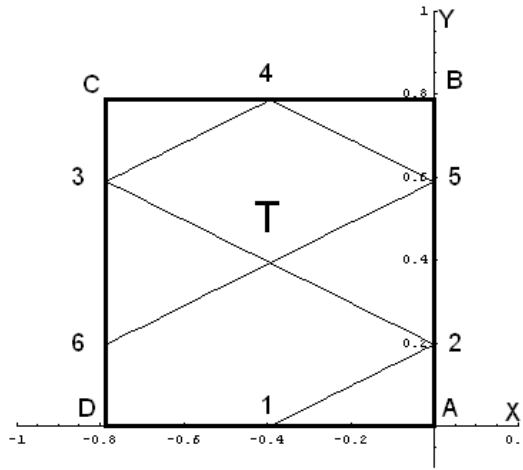

**Figure 1**. Trajectory of a light ray that is reflecting on the sides of a cuadrilateral. The refractive index n(x,y) is constant.

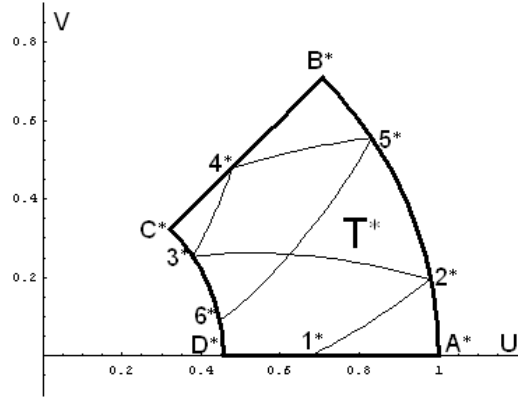

**Figure 2**. The image of Figure 1 by the conformal mapping $w = f(z) = ie^{-z}$. The refractive index that fills the interior of the curvilinear polygon is given by
$$n(u,v) = \frac{1}{\sqrt{u^2 + v^2}}.$$

## 3. Construction of an "optical torus" from the polygonal billiard

### 3.1 The mapping

The first task is to find a holomorphic function that maps the interior of our polygon, P, into the interior of a rectangle R. To that end we use two Schwartz-Christoffel (SC) mappings:

  **1.** f: Upper half plane → Polygon

Let´s suppose that we have in the z-plane a n-gon, P, with vertices $P_i$ at points $z_i$ (i=1,…,n) and interior angles at those vertices given by $\alpha_i \cdot \pi$. Suppose that the points $P_i$ are the image, under f, of the n points $A_i$, with coordinates $a_i$, situated in increasing order in the real axis of the η-plane (Fig.(3)). A transformation that maps the interior of the upper-half plane H of the η-plane onto the interior of the polygon P, and the real axis into the boundary of P is given by a holomorphic function, z = z (η), whose derivative is the SC formula [3]

$$\frac{dz}{d\eta} = A(\eta - a_1)^{\alpha_1 - 1}(\eta - a_2)^{\alpha_2 - 1} .. (\eta - a_n)^{\alpha_n - 1} \qquad (7)$$

In (7) we choose those branches of the functions $(\eta - a_i)^{\alpha_i - 1}$ which are direct analytic continuation of the real functions $(x - a_i)^{\alpha_i - 1}$ of the real variable $x > a_i$. So the function defined by (7) is a single-valued holomorphic function in H. The n constants $\alpha_i$ are subjected to the condition $\sum \alpha_i = n - 2$. Also we have that

$$\left|\frac{dz}{d\eta}\right|_{\eta=a_i} = \begin{cases} \infty & \text{if } \alpha_i < 1 \\ 0 & \text{if } \alpha_i > 1 \end{cases}$$

So the points $a_i$ lying on the real axis are singularities of the derivative of f. However it can be proved that f itself is continuous on this axis.

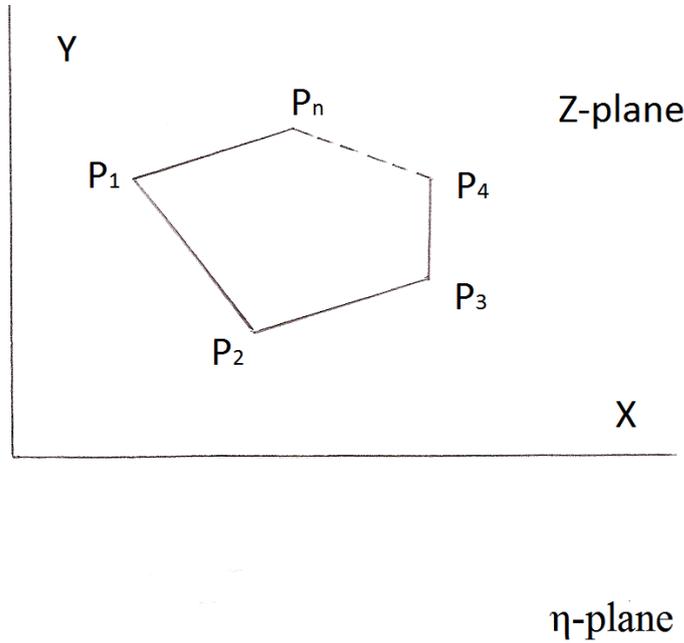

**Figure 3**. The SC transformation $z = z(\eta)$ given by Equation (7) transforms the upper-half $\eta$-plane into the interior of the polygon in the z-plane.

2. g: Upper half plane → Rectangle.

Applying the SC formula (7) to a rectangle, we obtain the transformation, g, which maps H onto the interior of a rectangle, R, in the w-plane (Fig. (4)). The mapping is given by

$$w = \int_0^\eta \frac{d\zeta}{\sqrt{(1-\zeta^2)(1-k^2\zeta^2)}} \equiv F(\eta\,;k) \quad (k < 1), \tag{8}$$

where $F(\eta\,;k)$ is the elliptic integral of the first type with parameter k. The real axis is transformed into the contour of the rectangle of Fig. (4) with the correspondences:

$$C_1(c_1 = -1/k) \leftrightarrow R_1(w_1 = -a - ib)$$
$$C_2(c_2 = -1) \leftrightarrow R_2(w_2 = -a)$$
$$C_3(c_3 = 1) \leftrightarrow R_3(w_3 = a)$$
$$C_4(c_4 = 1/k) \leftrightarrow R_4(w_4 = a + ib),$$

where a and b are given in terms of the complete elliptic integral of the first type K(k):

$$a = \int_0^1 \frac{dx}{\sqrt{(1-x^2)(1-k^2x^2)}} \equiv K(k) \quad , \quad b = \int_1^{1/k} \frac{dx}{\sqrt{(1-x^2)(1-k^2x^2)}} = K(k'), (k'^2 = 1-k^2)$$

The derivative of the mapping is singular at the four points $c_i$ (i = 1, 2, 3, 4) corresponding to the vertices of the rectangle:

$$\frac{dw}{d\eta} = (1-\eta^2)^{-1/2}(1-k^2\eta^2)^{-1/2} \qquad (9)$$

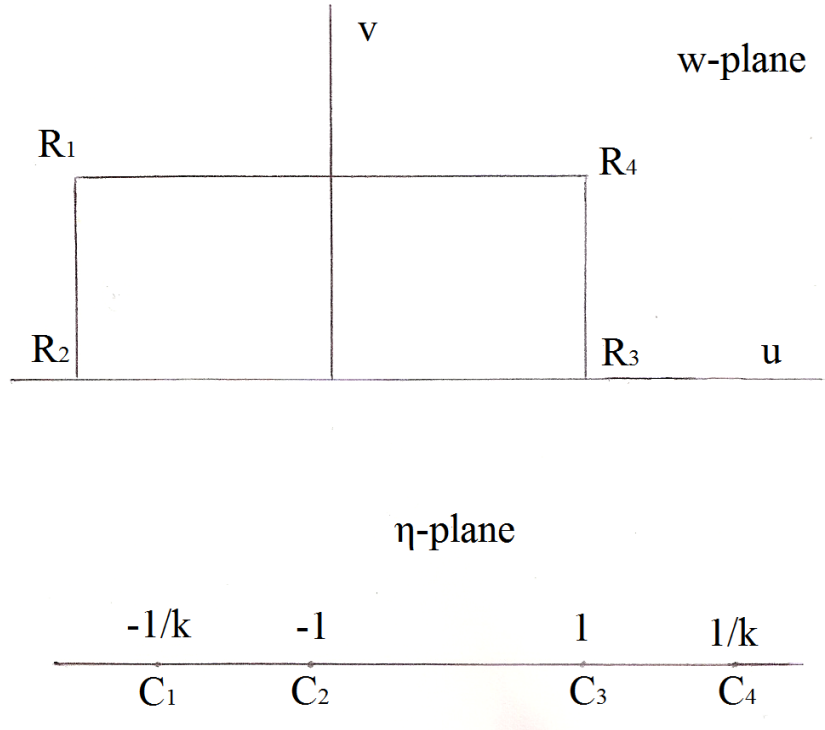

Figure 4. The SC transformation $w = w(\eta)$ given by Equation (8) transforms the upper-half $\eta$-plane into the interior of the rectangle in the w-plane.

### 3. h: Polygon → Rectangle

The holomorphic function, h, which transform the polygon P into the rectangle is $h = g \circ f^{-1}$. Now we take into account the fact that, given the polygon P, we can freely choose the position of three of the values of the $a_i$. The remaining n-3 values: $a_4, a_5, ..., a_n$ are the determined by the polygon. So we make the following choice:

$$a_1 = c_1 = -1/k, \; a_2 = c_2 = -1, \; a_3 = c_3 = 1 \qquad (10)$$

With this choice, h maps the three points $P_1, P_2, P_3$ of the polygon into the three ones $R_1, R_2, R_3$ of the rectangle. Moreover let $P_i$ be one of the remaining n-3 vertices: $P_4... P_n$. If $a_i < 1/k$ then $P_i$ is mapped into the interior of the vertical side $R_3$-$R_4$. On the contrary, if $a_i > 1/k$ then its image lies on the interior of the horizontal side $R_4$-$R_1$. So, given the polygon, we only have the freedom to choose a value for the parameter k, which determines the aspect ratio of the rectangle, and the vertex of the polygon ($P_1$) which is the image by f of $a_1 = -1/k$.

## 3.2 The refractive index

To obtain the refractive index, n(u,v), that must fill the rectangle R of the w-plane we observe that

$$n(w) = n(u,v) = \left|\frac{dz}{dw}\right| = \left|\frac{dz}{d\eta}\right|\left|\frac{d\eta}{dw}\right| \quad (11)$$

By substituting (7), (9) and (10) in (11):

$$n(w) = A\left|\eta + \frac{1}{k}\right|^{\alpha_1 - \frac{1}{2}} \cdot |\eta + 1|^{\alpha_2 - \frac{1}{2}} \cdot |\eta - 1|^{\alpha_3 - \frac{1}{2}} \cdot \left|\eta - \frac{1}{k}\right|^{\frac{1}{2}} \cdot \prod_{i=4}^{n} |\eta - a_i|^{\alpha_i - 1}, \quad (12)$$

where η can be expressed, according to Eq. (8), in terms of w by means of the inverse of the elliptic function of the first type (a Jacobian elliptic function):

$$\eta = F^{-1}(w;k) \equiv sn(w) \quad (13)$$

The function sn(w) has a simple pole at the middle point, $w_p = b \cdot i$, of the upper side of the rectangle R. This is due to the fact that g maps $\eta = \infty$ into $w_p$. Obviously, when the polygon is a triangle, the last n-3 products of Eq. (12) are absent and $\alpha_3 = 1 - (\alpha_1 + \alpha_2)$.

Equation (12) clearly shows, as previously observed, that the logarithm of n is a harmonic function. In terms of the variable η, this logarithm is a linear combination of the logarithms of the distances from η to the singularities $a_i$ lying in the boundary of the upper-half plane. If P is a rational polygon then all of the coefficients of that linear combination are also rational.

## 3.3 The "optical" tessellation

Let´s imagine that we manufacture a thin rectangular optical sheet whose refractive index is given by (12) and (13). Moreover, suppose that the boundaries of P and R are perfects mirrors. Then, to each rectilinear polygonal light trajectory in P, whose refractive index is some constant, there corresponds a curvilinear one in the rectangular sheet R, and vice versa. This correspondence is given by the holomorphic function h. The refractive index of R is singular at the following points: the four corners; the n-3 points lying on the interior of the two consecutives sides, $R_3$-$R_4$ and $R_4$-$R_1$, and at the pole, D, of sn(w) situated in the middle of the upper side $R_4$-$R_1$.

The next step is to unfold the trajectories in R by the known method [4] of reflecting R with respect to the side on which the light ray is reflecting. In Fig.(5) we have unfolded a certain fragment of the image in R of a hypothetical trajectory in a hexagon. The points A, B, and C are the images, on the boundary of R, of the points $P_4$, $P_5$ and $P_6$, respectively. The fundamental domain, F, of our "optical" tessellation is formed by R and its three "reflected" copies, that have been denoted by $R_r$ (R-right), $R_u$(R-up) and $R_{r-u}$(R-right-up). This operation is equivalent to extending to F the metric defined on R by means of successive reflections on its boundary.

Let $n_e(u,v)$ be the extended conformal factor in F of the one n(u,v) defined in R by (12) and (13). That is to say:

$$n_e(u,v) = \begin{cases} n(2a-u, v) & \text{if } (u,v) \in R_r \\ n(u, 2b-v) & \text{if } (u,v) \in R_u \\ n(2a-u, 2b-v) & \text{if } (u,v) \in R_{r-u} \end{cases} \qquad (14)$$

We know that n(u,v) is continuous in R and its boundary (with the exception of the singular points), so the extension $n_e(u,v)$ is also continuous in F with the exception of the singularities situated in the boundary of R and its reflected copies (see Fig.(5)).

The physical manifestation of the continuity of the extended refractive index $n_e(u,v)$ is the smoothness of the unfolded trajectory. This is due to the fact that the law of reflection is valid in R as well.

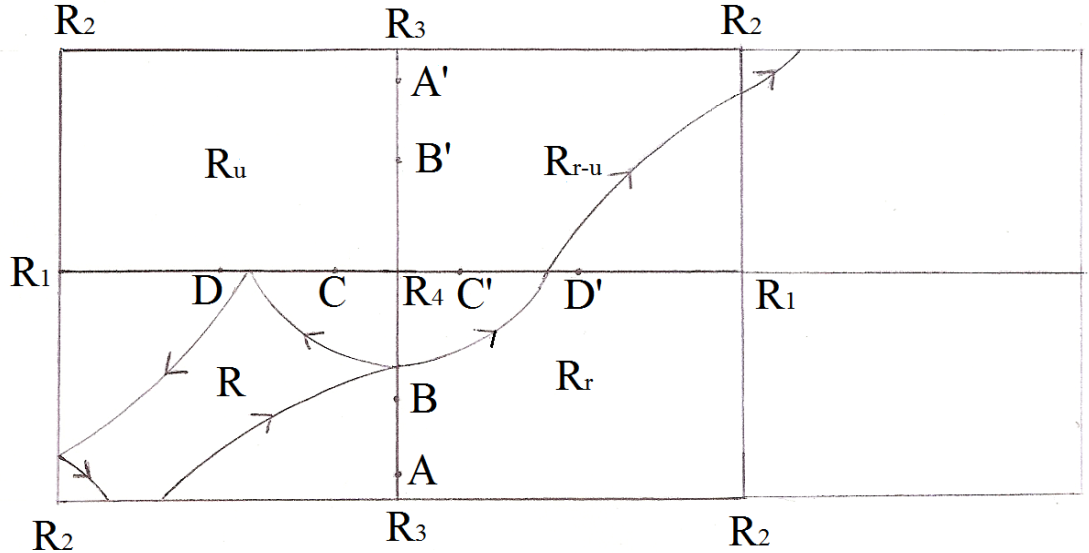

**Figure 5**. Unfolding the curvilinear trajectory inside R. The fundamental domain, F, of the "optical" tessellation is formed by the four rectangles : R, $R_r$, $R_u$ and $R_{r-u}$.

### 3.4 The "optical" torus

Let´s identify, in the standard way, the two pairs of parallel sides of the fundamental domain F. The result is a torus, T, formed by welding of four copies of our optical sheet (Fig. (6)). This "optical torus" has some "impurities". They are the points where the refractive index is singular: the four points $P_1$ to $P_4$ and the eight points A, A', B, B', C, C' and D, D' corresponding to the images of the vertices $P_4$, $P_5$ and $P_6$ of the hexagon and the pole of sn (w), respectively. In the general case of a n-gon, the number of singular points is 4+2+2 (n-3) =2n. They are situated on two pairs of opposite meridians and parallels of T, on which n(u,v) is continuous but not differentiable. So T is a one-holed 2n punctured torus endowed with the flat metric

$$ds^2 = n_e(u,v)^2 (du^2 + dv^2),$$

where the conformal factor $n_e(u,v)$ is given by Eqs. (12), (13) and (14). We could say that, in a way, all the complexity of the billiards trajectories in the polygon is "encrypt-ed" in T.

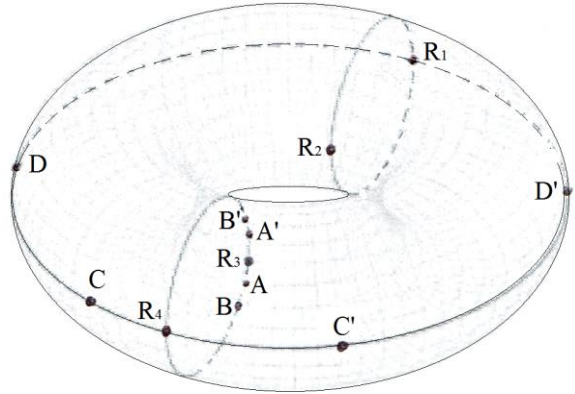

**Figure 6.** The torus of the optical tessellation (Fig. (5)) that correspond to a certain hexagonal billiard. There are 2*6=12 singular points. The location of the six points $R_1$, $R_2$, $R_3$, $R_4$, D and D' are independent of the number of sides of the polygon. The remaining 2n-6 points (cero for a triangle) are situated on the parallel $R_1$-$R_4$ and the meridian $R_3$-$R_4$.

### 3.5 A periodic potential

There is a useful geometrical optics formalism [5] in which the equation that governs the trajectories of light rays assumes the form of the Newton's second law:

$$\frac{d^2 \vec{x}}{d a^2} = \nabla(\frac{1}{2}n^2) \qquad (15)$$

In the previous equation $\vec{x}$ is the position of the light ray that moves in the plane $R^2$ (equipped with the standard metric $ds^2 = du^2 + dv^2$) through a region of varying refractive index $n(\vec{x})$. The (stepping) parameter, a, is defined by $|d\vec{x}/da| = n$. This parameter is the optical analog of time t. Also, the optical analog of the potential energy U is

$$U(\vec{x}) = -\frac{1}{2}n(\vec{x})^2 \qquad (16)$$

If the particle has a unit mass, its "kinetic energy" is $T = \frac{1}{2}|d\vec{x}/da|^2$ and the "total energy" is cero:

$$E = T + U = \frac{1}{2}n^2 + (-\frac{1}{2}n^2) = 0 \qquad (17)$$

So we have encountered the problem of studying the unbounded motion of a particle with cero total energy in a bidimensional periodic potential that is always negative (it can be null in certain isolated points). The singular points at which n(w) is infinity ( see Eq. (12)) corresponds to the infinite negative potential wells of $U(\vec{x})$.

Alternatively, Eq. (15) describes the movement of a Newtonian particle in a torus on which the potential energy field is given by (16).

## References

[1] M. Born, E. Wolf. *Principles of Optics*. Pergamon Press, 1970.

[2] B. Doubrovine, S. Novikov, A. Fomenko. *Géométrie contemporaine. Méthodes et applications*. Premièr Partie. Éditions Mir. Moscou, 1985.

[3] A. Sveshnikov, A. Tikhonov. *The Theory of functions of a complex variable*. Mir Publishers, 1978.

[4] S. Tabachnikov. *Geometry and Billiards*. Student Mathematical Library. Volume 30. AMS. 2005.

[5] J. Evans. *Simple forms of rays in gradient-index lenses*. Am. J. Phys **58** (8), August 1990.